\documentclass[11pt]{article}
\usepackage{graphicx}

\newfont{\bb}{msbm10 at 11pt}
\def\r{\hbox{\bb R}}
\newcommand{\T}{\hbox{\bf t}}
\newcommand{\G}{\hbox{\bf G}}
\newcommand{\N}{\hbox{\bf n}}
\newcommand{\B}{\hbox{\bf b}}
\newcommand{\C}{\hbox{\bf c}}
\newcommand{\E}{\hbox{\bf e}}
\newcommand{\X}{\hbox{\bf X}}
\newcommand{\e}{\hbox{\bf E}}
\newcommand{\h}{\hbox{\bf H}}

\newtheorem{theorem}{Theorem}[section]
\newtheorem{corollary}{Corollary}
\newtheorem{definition}{Definition}[section]
\newtheorem{remark}{Remark}[section]
\begin{document}

\title{Linear Weingarten surfaces foliated by circles\\ in Minkowski space }
\author{\"{O}zg\"{u}r Boyacioglu Kalkan\\ Mathematics Department\\ Afyon Kocatepe Universtiy\\Afyon 03200 Turkey\\
email: 	bozgur@aku.edu.tr\\
\linebreak \\
Rafael L\'opez\footnote{Partially supported by MEC-FEDER
 grant no. MTM2007-61775 and
Junta de Andaluc\'{\i}a grant no. P06-FQM-01642.}\\Departamento de Geometr\'{\i}a y Topolog\'{\i}a\\ Universidad de Granada\\ 18071 Granada, Spain\\
email: rcamino@ugr.es\\
\linebreak \\
Derya Saglam\\Mathematics Department\\ Afyon Kocatepe Universtiy\\Afyon 03200 Turkey}

\date{}
\maketitle


\begin{abstract} In this work,
we study spacelike surfaces in Minkowski space $\e_1^3$  foliated by  pieces of circles and
that satisfy a linear Weingarten condition of type
$a H+b K=c$, where $a,b$ and $c$ are constant and  $H$ and $K$ denote
the mean curvature and the Gauss curvature respectively. We show that such surfaces must be surfaces of revolution or surfaces with constant mean curvature $H=0$ or  surfaces with constant Gauss curvature $K=0$.
\end{abstract}

%%%%%%%%%%%%%%%%%%%%%%%%%%%%%%%%%%%%%%%%%%%%%
\section{Introduction and results}\label{intro}
%%%%%%%%%%%%%%%%%%%%%%%%%%%%%%%%%%%%%%%%%%%%

Let $\e_1^3$ be the Minkowski three-dimensional space, that is, the real vector space $\r^3$ endowed with
the scalar product
$$\langle,\rangle=(dx_1)^2+(dx_2)^2-(dx_3)^2,$$
where $(x_1,x_2,x_3)$ denote the usual coordinates in $\r^3$. An immersion $x:M\rightarrow\e_1^3$ of a surface $M$ is  called {\it
spacelike} if the induced metric $x^{*}\langle,\rangle$ on $M$  is
 a Riemannian metric. In this
paper, we study spacelike surfaces that satisfy a relation of type
\begin{equation}\label{w1}
a H+b K=c,
\end{equation}
where $H$ and $K$ are the mean curvature and  the Gauss curvature of
$M$ respectively,  and $a, b$ and $c$ are constant with
$a^2+b^2\not=0$.  We say that $M$ is a {\it linear Weingarten
surface}.  These surfaces generalize the surfaces with constant mean curvature ($b=0$ in (\ref{w1})) and the surfaces with constant Gauss curvature
($a=0$ in (\ref{w1})).   In Euclidean space, there is a great amount of literature of Weingarten surfaces, beginning with  works of Chern, Hartman, Winter and Hopf in the fifties, and more recently in \cite{gmm,ks,lo2,lo3,rs}, without being a complete bibliography.

In order to look for  examples of linear Weingarten
spacelike surfaces in $\e_1^3$, it  is natural to assume some hypothesis about
the geometry of the surface. A simple condition is that the surface
is rotational. In such case, Equation (\ref{w1}) is an ODE of second
order given in terms of the generating curve of $M$. In a more general scene, we consider surfaces given by a foliation of circles. Following terminology due to Enneper, we give the next definition.

\begin{definition} A cyclic surface in Minkowski space $\e_1^3$ is a surface determined
 by a smooth uniparametric family of  circles.
\end{definition}

We also say that the surface is foliated by circles.
As in Euclidean space, by a circle in $\e_1^3$ we mean a planar curve with
constant curvature. In particular, given a cyclic surface there exists a uniparametric
family of planes of $\e_1^3$ whose intersection with $M$ is a circle. Since the
circles are contained in a spacelike surface, each circle of the
foliation must be a spacelike curve. However, the planes containing the circles can be of any causal type.

Our work is motivated by the following fact. In Minkowski space $\e_1^3$ there are cyclic spacelike surfaces with $H=0$ (or $K=0$) that are not rotational surfaces. For the maximal case ($H=0$) these surfaces are foliated by circles in parallel planes and they represent in Minkowski  ambient the same role as the classical Riemann examples of minimal surfaces in Euclidean space.  These surfaces appeared for the first time in the literature in  \cite{lo1} and they have been origin of an extensive study in recent years: \cite{al,fl,fls,ky,lls}. In the same sense,  non-rotational cyclic surfaces with constant Gauss curvature $K=0$ appeared in \cite{lo11}. See also Remark \ref{r-1}.  Thus, it is natural to ask if besides these examples, there exist other cyclic surfaces in the family of linear Weingarten  surfaces of $\e_1^3$. If we compare with what happens in Euclidean space, the difficulty in $\e_1^3$ is the variety of possible cases that can appear since the plane containing the circle can be of spacelike, timelike or lightlike type. Because we are looking for new cyclic linear Weingarten surfaces, we will exclude throughout this work the known examples corresponding to $H=0$ ($b=c=0$) and $K=0$ ($a=c=0$).

In the case that the planes of the foliation are parallel, we prove:

\begin{theorem}\label{t1}  Let $M$ be a spacelike cyclic surface in $\e_1^3$ and we assume that the circles of the foliation lie in parallel planes. If   $M$ is a linear Weingarten surface, then $M$ is a surface of revolution, or $H=0$ or $K=0$.
\end{theorem}

In Minkowski space $\e_1^3$ there are spacelike surfaces that play the same role as spheres in Euclidean space. These surfaces are   the pseudohyperbolic surfaces.
After an isometry of $\e_1^3$, a pseudohyperbolic surface of radius $r>0$ and centered at $x_0\in\e_1^3$ is  parametrized as
$$\h^{2,1}(r,x_0)=\{x\in\e_1^3;\langle x-x_0,x-x_0\rangle=-r^2\}.$$
From the Euclidean viewpoint, $\h^{2,1}(r,O)$ is the hyperboloid of two sheets $x_1^2+x_2^2-x_3^2=-r^2$ which is obtained by rotating the hyperbola $\{x_1^2-x_3^2=r^2,x_2=0\}$ with respect to the $x_3$-axis. This surface is spacelike with constant mean curvature $H=1/r$ and with constant Gauss curvature $K=1/r^2$. In particular, $\h^{2,1}(r,x_0)$ is a linear Weingarten surface: exactly, there are many choices of constants $a,b$ and $c$ that satisfy (\ref{w1}).  Although this surface is rotational, any uniparametric family of (non-parallel) planes intersects $\h^{2,1}(r,x_0)$ in circles. Taking account this fact about the pseudohyperbolic surfaces, our next result establishes:

\begin{theorem}\label{t2}  Let $M$ be a spacelike cyclic surface in $\e_1^3$. If $M$ is a linear Weingarten surface,  then $M$ is a pseudohyperbolic surface or the planes of the foliation are parallel.
\end{theorem}

As consequence of the above two results, we have

\begin{corollary} The only non-rotational spacelike cyclic surfaces that are linear Weingarten surfaces are the Riemann examples of maximal surfaces \cite{lo1} and a family of surfaces with $K=0$ described in \cite{lo11}.
\end{corollary}

The proof of Theorems \ref{t1} and \ref{t2} involves   long and complicated algebraic computations that have been possible check them by using a symbolic program such as Mathematica.

Finally, we point out that Theorems \ref{t1} and \ref{t2}  hold for linear Weingarten cyclic {\it timelike} surfaces of $\e_1^3$. The proofs are similar and we do not included them in the present paper, althout they can easily carried. We only remark the two differences with the spacelike case. First, it appears a new parametrization of circle, which it is a {\it timelike} curve.  On the other hand, the first fundamental form in classical notation $W=EG-F^2$ is negative, in contrast to the spacelike case, that it is positive. However, the key fact that we use in our proofs is that the metric is non-degenerate, that is, $W\not=0$ on the surface, independent if $W>0$ (spacelike) or $W<0$ (timelike).

%%%%%%%%%%%%%%%%%%%%%%%%%%%%%%%%%%%%%%%%%%%%%%%
\section{Preliminaries}
%%%%%%%%%%%%%%%%%%%%%%%%%%%%%%%%%%%%%%%%%%%%%%%%

A vector $v\in\e_1^3$ is said spacelike if $\langle v,v\rangle>0$ or
$v=0$, timelike if $\langle v,v\rangle<0$ and lightlike if $\langle
v,v\rangle=0$ and $v\not=0$.  A plane $P\subset\e_1^3$ is said
spacelike, timelike or lightlike if the induced metric on $P$ is a
Riemannian metric (positive definite), a Lorentzian metric (a metric
of index $1$) or a  degenerated metric, respectively. This is
equivalent that any orthogonal vector to $P$ is timelike, spacelike
or lightlike respectively.

 Consider
$\alpha:I\subset\r\rightarrow\e_1^3$ a parametrized regular
  curve in $\e_1^3$. We say that $\alpha$
is spacelike  if $\alpha'(t)$ is  spacelike for all $t\in I$. We can
reparametrize $\alpha$ by a parameter $s$ such that
$\langle\alpha'(s),\alpha'(s)\rangle=1$ for any $s\in I$. Then one
can define a Frenet trihedron at each point. The differentiation of
the Frenet frame allows to define the curvature $\kappa$ and the
torsion $\tau$ of $\alpha$. See \cite{ku,on}. Motivated by what
happens in Euclidean ambient, we give the following definition:

\begin{definition} A spacelike curve in Minkowski space $\e_1^3$ is a planar curve with constant curvature.
\end{definition}

We describe the classification of spacelike circles in $\e_1^3$. This classification depends on the causal character of the plane $P$ containing the circle.  After
an isometry of the ambient space $\e_1^3$,  a circle parametrizes as
follows:
\begin{enumerate}
\item If $P$ is the horizonal plane $x_3=0$, the circle is given by
$$\alpha(s)=r\Big(\cos(\frac{s}{r}),\sin(\frac{s}{r}),0\Big),\ \ r>0.$$
In this case, the curve is a Euclidean horizontal circle.
\item If $P$ is the vertical plane $x_1=0$, then
$$\alpha(s)=r\Big(0,\sinh(\frac{s}{r}),\cosh(\frac{s}{r})\Big),\ \ r>0.$$
The curve describes  a hyperbola in a vertical plane.
\item If $P$ is the  plane $x_2-x_3=0$, then
$$\alpha(s)=\Big(s,r\frac{s^2}{2},r\frac{s^2}{2}\Big), \ \ r>0.$$
The curve is a parabola in $P$.
\end{enumerate}

A surface $M$ in $\e_1^3$ is a surface of revolution (or rotational
surface) if there exists a straight line $l$ such that $M$ is
invariant by the rotations that leave $l$ pointwise fixed. In particular, a rotational surface in $\e_1^3$ is formed by a uniparametric family
of circles of $\e_1^3$ in parallel planes.

 Let $M$ be a spacelike
surface in $\e_1^3$. The spacelike condition is equivalent that any
unit normal vector $\G$ to $M$ is always timelike. Since any two
timelike vectors in $\e_1^3$ can not be  orthogonal, $\langle \G,(0,0,1)\rangle\neq 0$ on $M$. This shows that
$M$ is an orientable surface. As in Euclidean space, one can define
the mean curvature $H$ and the Gauss curvature of $K$ of $M$ as:
$$H=\frac12\ \mbox{trace}(d\G),\hspace*{1cm}K=\mbox{det}\ (-d\G).$$
If we locally write the immersion as  $\X(u,v)$, with $(u,v)$ in
some planar domain, then the following formulae are well-known
\cite{we}:
$$H=\frac12\ \frac{eG-2f F+gE}{EG-F^2},\hspace*{1cm}K=\frac{e\ g-f^2}{EG-F^2},$$
where $\{E,F,G\}$ and $\{e,f,g\}$ are the coefficients of the first and second fundamental
 forms respectively of the immersion according to the orientation $$\G=\frac{\X_u\wedge \X_v}{|\X_u\wedge \X_v|}.$$
 Here $\wedge$ stands for the Lorentzian cross product and
$$E=\langle \X_u,\X_u\rangle,\quad F=\langle \X_u,\X_v\rangle,\quad
G=\langle \X_v,\X_v\rangle.$$
$$e=\langle \G, \X_{uu}\rangle,\quad f=\langle\G,\X_{uv}\rangle,
\quad g=\langle \G,\X_{vv}\rangle.$$
Denote $W=:EG-F^2=|\X_u\wedge \X_v|^2$. This function is positive because the immersion is
spacelike. From the expressions of $H$ and $K$, we have
$$G[\X_u,\X_v,\X_{uu}]-2F [\X_u,\X_v,\X_{uv}]+E[\X_u,\X_v,\X_{vv}]=2H W^{3/2}$$
$$[\X_u,\X_v,\X_{uu}] [\X_u,\X_v,\X_{vv}]-[\X_u,\X_v,\X_{uv}]^2= K W^2$$
 and $[,,]$ denotes the determinant of three vectors: $[v_1,v_2,v_3]=\mbox{det}(v_1,v_2,v_3)$.

%%%%%%%%%%%%%%%%%%%%%%%%%%%%%%%%%%%%%%%%%%%%%%%%%%%

\section{Proof of Theorem \ref{t1}}

%%%%%%%%%%%%%%%%%%%%%%%%%%%%%%%%%%%%%%%%%%%%%%%%%%%%

We consider a spacelike surface $M\subset\e_1^3$ parametrized by circles in parallel planes. We distinguish three cases according to the causal character of the planes of the foliation.

\subsection{The planes are spacelike}

After a rigid motion in $\e_1^{3},$ we may assume the planes are
parallel to the plane $x_{3}=0$.  The circles are horizontal
Euclidean circles and $M$ can be parametrized by
\begin{equation}\label{para1}
\X(u,v)=(f(u),g(u),u)+r(u)(\cos {v},\sin {v},0),
\end{equation}
where $f,g,r>0$ are smooth functions in some $u$-interval $I$.  With this parametrization, $M$ is a surface of revolution
if and only if $f$ and $g$ are constant functions.

The Weingarten relation $aH+bK=c$ writes as:
\begin{eqnarray}\label{31}
&&a\ \frac{G[\X_u,\X_v,\X_{uu}]-2F [\X_u,\X_v,\X_{uv}]+E[\X_u,\X_v,\X_{vv}]}{2W^{3/2}}\\
& & +b\ \frac{[\X_u,\X_v,\X_{uu}]
[\X_u,\X_v,\X_{vv}]-[\X_u,\X_v,\X_{uv}]^2}{W^2}=c.
\end{eqnarray}
%%%%%%%%%%%%%%%%%%%%%%%%%%%

\subsubsection{Case $c=0$.}

%%%%%%%%%%%%%%%%%%%%%%%%%%%
Equation (\ref{31}) writes as
\begin{eqnarray*}
& &a^2\ \Big(G[\X_u,\X_v,\X_{uu}]-2F [\X_u,\X_v,\X_{uv}]+E[\X_u,\X_v,\X_{vv}]\Big)^2W^2\\
& & -4b^2\  \Big([\X_u,\X_v,\X_{uu}]
[\X_u,\X_v,\X_{vv}]-[\X_u,\X_v,\X_{uv}]^2\Big)^2=0.
\end{eqnarray*}
 Without loss of generality,
we assume $4b^2=1$. If we compute the above expression with the parametrization $\X(u,v)$,  we
obtain an expression
\begin{equation}
\sum_{j=0}^{4}A_{j}(u)\cos {(jv)}+B_{j}(u)\sin {(jv)}=0.  \label{p0}
\end{equation}%
Then the functions $A_{j}$ and $B_{j}$ on $u$ must vanish on $I$. By contradiction, we assume that $M$ is not rotational. Then $f'$ or $g'$ does not vanish in some interval.

\textbf{1.} We consider the cases that one of the functions $f$ or
$g$ is constant. For simplicity we  consider $f'=0$ in
some interval. Then $g'\not=0$. The coefficient $A_{4}$ writes as
$$A_{4}=\frac{1}{8}a^2r^{6}g'^2(rg''-2r'g')^{2}.$$
As $g'\neq 0,$ we have that $rg''-2r'g'=0$. Then $g'=\lambda r^{2}$ for some positive
constant $\lambda \neq 0.$ Now
\begin{eqnarray*}
A_{2} &=&\frac{1}{2}\lambda^{2}r^{8}(4r'^2-a^2r^{2}A^{2}) \\
B_{1} &=&2\lambda r^{7}r'(a^2rA^{2}-2r'')
\end{eqnarray*}%
where $$ A=-1+\lambda ^{2}r^{4}+r'^2-rr''.$$
Firstly, from Equation $B_1=0$, we discard the case that $r$ is a constant. In such case, the coefficient $E$ of the first fundamental form  vanishes. As a consequence and since $r'\not=0$, the combination of $A_{2}=0$ and $B_{1}=0$ leads to
that function $r$ satisfies the ordinary differential equation
$2r'^2-rr''=0.$ Then
$$ r(u)=\frac{c_{2}}{u+c_{1}}, \ \
c_{1},c_{2}\in \r.$$
Now $A_{2}=0$ gives a polynomial equation on $u$
given by
$$-4(u+c_{1})^{6}+a^2((u+c_{1})^{4}+c_{2}^{2}-\lambda
^{2}c_{2}^{4})^{2}=0.$$
 In particular, the leading coefficient $a^2$ must vanish:  contradiction. This means that
the assumption that $f$ is constant is impossible.

\textbf{2.} We assume that both $f$ and $g$ are not constant
functions. Then $f',g'\neq 0$. The coefficient $B_{4}$ yields:
$$(-4f'g'r'+rf'f''+rf'g'')(-2f'^2r'+2g'^2r'+rf'f''-r f'g'')=0.$$
We distinguish     two cases:
\begin{enumerate}
\item Assume $-4f'g'r'+rg'f''+rf'g''=0$. Then
$$f''=\frac{4f'g'r'-rf'g''}{rg'}.$$
Now $A_{4}=0$ gives
$$a^2r^{6}(f'^2+g'^2)^{2}(-2g'r'+rg'')^2=0,$$
 that is,
\begin{equation}\label{eq-g2}
g''=\frac{2g'r'}{r}.
\end{equation}
This implies $g'=\lambda r^{2}$ with $\lambda >0$. Analogously,
$f'=\mu r^{2}$, $\mu >0$. The computation of $B_{2}$ and $B_{1}$
leads to
$$B_2=\lambda\mu r^8 (a^2r^2 A^2-4r'^2),$$
$$B_1=2\lambda r^7 r'(a^2rA^2-2r''),$$
where the value of $A$ is
$$A=-1+(\lambda ^{2}+\mu ^{2})r^{4}+r'^2-r''.$$
Equation $B_1=0$ gives the possibility  $r'=0$, that is, $r$ is a constant function. In such case, $B_2=\lambda\mu a^2(-1+(\lambda^2+\mu^2)r^4)^2$. The computation of the coefficient  $E$ of the first fundamental form gives $E=0$: contradiction. Thus, we can assume that $r'\not=0$.

By combining $B_2=B_1=0$, we obtain $r r''=2 r'^2$. Solving this equation, we have
$$r(u)=\frac{c_2}{u+c_1},\ \ c_1,c_2\in\r.$$
The coefficient $B_2$  writes now as as polynomial on $u$ and from $B_2=0$ we conclude
$$a^2(u+c_2)^8-4(u+c_2)^6+2a^2c_2^4(1-c_2^2(\lambda^2+\mu^2))(u+c_2)^4+a^2c_2^4(1-c_2^2(\lambda^2+\mu^2))^2=0.$$
The leading coefficient must vanish, that is, $a^2=0$: contradiction.

\item Assume  $-2f'^2r'+2g'^2r'+rf'f''-rg'g''=0$.  From here, we
obtain $f''$ and putting it into $A_{4}$, it gives
$$A_{4}=-\frac{a^2r^6(f'^2+g'^2)^2}{8f'^2}
(rg''-2g'r')^2.$$
 Then $rg''-2g'r'=0$ and we now are  in the position of the above case  (\ref{eq-g2}) and this finishes the proof.
\end{enumerate}

%%%%%%%%%%%%%%%%%%%%%%%%%%%%%%5
\subsubsection{Case $c\neq 0$}
%%%%%%%%%%%%%%%%%%%%%%%%%%%%

The computation of $A_{8}$ and $B_{8}$ gives respectively:
\begin{eqnarray*}
A_{8} &=&-\frac{1}{32}c^{2}r^{8}(f'^8-28f'^6 g'^2+70f'^2g'^6+g'^8) \\
B_{8} &=&\frac{1}{4}c^{2}r^{8}f' g'(-f'^6-7f'^4g'^2-7f'^2 g'^4+g'^6)
\end{eqnarray*}
Since $\alpha (u)=(f(u),g(u),0)$ is not a constant planar curve, we
parametrize it by the arc-length, that is, $(f(u),g(u))=(x(\phi
(u),y(\phi (u)),$ where
$$f'(u)=\phi'(u)\cos (\phi (u)),\ \ g'(u)=\phi'(u)\sin (\phi (u)),\ \
\phi'^2=f'^2+g'^2.$$ With this change of variable, the functions
$A_{8}$ and $B_{8}$ write now as:
\begin{eqnarray*}
A_{8} &=&-\frac{1}{32}c^{2}r^{8}\phi'^8\cos (8\phi (u)). \\
B_{8} &=&-\frac{1}{32}c^{2}r^{8}\phi'^8\sin (8\phi (u)).
\end{eqnarray*}
As $c\neq 0$ and $r>0$, we conclude that $\phi'=0$ on some interval. Therefore $f'^2+g'^2=0$, which means that
$\alpha $ is a constant curve, obtaining a contradiction.
This finishes the proof of Theorem \ref{t1} for the case that the
planes are spacelike.

%%%%%%%%%%%%%%%%%%%%%%%%%%%%%%%%%%%%%%%%%%%%%
\subsection{The planes are timelike }
%%%%%%%%%%%%%%%%%%%%%%%%%%%%%%%%%%%%%%%%%%
Let $M$ be a  linear Weingarten spacelike surface foliated by pieces
of circles in parallel timelike planes. After a motion in $\e_1^{3},$
we can suppose that these planes are parallel to the plane
$x_{1}=0.$ In this case we  parametrize the surface by
\begin{equation}\label{para2}
\hbox{\bf X}(u,v)=(u,f(u),g(u))+r(u)(0,\sinh v,\cosh v),
\end{equation}
where $r>0, f$ and $g$  are smooth functions. This means that $M$ is formed by a uniparametric family of vertical hyperbolas. In order to conclude
that $M$ is rotational it suffices to prove that $f$ and $g$ are
constant.

%%%%%%%%%%%%%%%%%%%%%%%%%%%%%%%%%%%%%%%%%%%%%

\subsubsection{Case $c=0$.}
%%%%%%%%%%%%%%%%%%%%%%%%%%%%%%%%%%%%%%%%%%%%%

As in the case of spacelike planes, the reasoning is by contradiction. We assume that $f$ or $g$ are not constant, that is, the functions $f'$ or $g'$ do not vanish.

\textbf{1.} Firstly, we consider the cases that one of the functions
$f$ or $g$ is constant. For simplicity, we shall consider $f'=0$ in some interval. Then $A_{4}$ writes as
$$A_{4}=-\frac{1}{8}a^2r^{6}g'^2(-2r'g'+rg'')^2.$$
 As $g'\neq 0,$ we have that
$rg''-2r'g'=0.$ Then $g'=\mu r^{2}$ for some positive constant $\mu \neq 0.$ Now
$$A_{2}=-\frac{1}{2}\mu ^{2}r^{8}(4r'^2+a^2r^{2}A^{2}),$$
$$A_{1}=-2\mu r^{7}r'(2r''+a^2rA^{2}),$$
 where
$$A=-1+\mu ^{2}r^{4}-r'^2+rr''.$$
As $a^2>0$, Equation $A_{2}=0$ implies that $r$ is a constant function and $A=-1+\mu^2 r^4=0$. Then the computation of the coefficient $E$ of the first fundamental form yields $E=0$: contradiction.

\textbf{2.} We assume that both $f$ and $g$ are not constant
functions. Then $f', g'\neq 0$. The coefficient
$B_{4}$ yields:
$$(-4f'g'r'+rg'f''+rf'g'')(-2f'^2r'-2g'^2r'+rf'f''+rg'g'')=0.$$
We consider two cases.

\begin{enumerate}
\item If $-4f'g'r'+rg'f''+rf'g''=0$, then

\begin{equation}
f''=\frac{f'(4f'g'r'-rg'')}{rb'}.  \label{1}
\end{equation}
Now $A_{4}=0$ gives
$$A_{4}=\frac{a^2r^{6}(f'^2-g'^2)^{2}(-2g' r'+rg'')^{2}}{8g'^2}.$$

\begin{enumerate}
\item[(a1)]  If $f'^2-g'^2=0$ then $g'=\pm f'$. Let $g'=f'$ (the case $g'=-f'$ is similar). Then $g=f+c_1,$ $c_1\in \r$. Putting it into $A_{4}$ and $B_{4},$ we obtain
\begin{eqnarray*}
A_{4} &=&-a^2r^{6}f'^2(-2f'r'+rf'') \\
B_{4} &=&a^2r^{6}f'^2(-2f'r'+rf'')^{2}
\end{eqnarray*}
As $f'\not= 0$, then $2f'r'=rf''$. Then $f'=\lambda r^2$ for some positive constant $\lambda$. The computation of $A_2$ and $B_1$ give
\begin{eqnarray*}A_{2} &=&\lambda^{2}r^{8}(4r'^2+a^2r^{2}A^{2}) \\
B_{1} &=& 2\lambda r^{7}r'(2r''+a^2rA^{2})
\end{eqnarray*}
where $A=1-r'^2+rr''$.
From Equation $A_{2}=0$ and the value of $A,$ we discard the case that $r$
is constant function. The combination of $A_{2}=0$ and $B_{1}=0$ implies that the function $r$ satisfies  $2r'^2-rr''=0$. Then
$$r(u)=\frac{c_{2}}{u+c_{1}},\ \ c_{1},c_{2}\in \r.$$
But then $A_{2}=0$ gives a polynomial on $u$ given by
$$4(u+c_{1})^{6}+a^2((u+c_{1})^{4}+c_{2}^{2})^{2}=0$$
whose leading coefficient is $a^2$: contradiction.

\item[(a2)] If $rg''=2g'r'$ then $g'=\mu r^{2}$ with $\mu >0$. Using (\ref{1}), the
same occurs for $f$: $f'=\lambda r^{2}$, $\lambda >0$.  The computation
of $A_{2}$ and $A_{1}$ leads to
\begin{eqnarray*}
A_{2} &=&-\frac{1}{2}(\lambda ^{2}+\mu
^{2})r^{8}(a^2r^{2}A^{2}+4r'^2)
\\
A_{1} &=&2\lambda r^{7}r'(2r''+a^2rA^{2})
\end{eqnarray*}%
where the value of $A$ is now
$$A=-1+(-\lambda ^{2}+\mu ^{2})r^{4}+r'^2-rr''.$$
Equation $A_2=0$ implies that $r$ is a constant function and $(\lambda^2-\mu^2)r^4=-1$. This gives $E=0$: contradiction.
\end{enumerate}

\item  $-2f'^2r'-2g'^2r'+rf'f''+rg'g''=0$. We obtain $f''$ and putting it into $A_{4},$ we obtain
$$A_{4}=\frac{a^2r^{6}(f'^2-g'^2)^{2}(-2g'r'+rg'')^2}{8g'^2}.$$

\begin{enumerate}
\item[(b1)] If $f'^2-g'^2=0$ then $g'=\pm f'$.  Now we are in the position of the above case (a1).

\item[(b2)]If $rg''=2g' r'$ then $g'=\lambda r^{2}$ with $\lambda >0.$ Now we are in the
position of the above case (a2).
\end{enumerate}
\end{enumerate}
%%%%%%%%%%%%%%%%%%%%%%%%%%%%%%%%%5
\subsubsection{Case $c\neq 0$}
%%%%%%%%%%%%%%%%%%%%%%%%%%%%%%%%%%%%%%%%%

The computations of $A_{8}$ and $B_{8}$ give respectively:
\begin{eqnarray*}
A_{8} &=&-\frac{1}{32}c^{2}r^{8}(f'^8+28f'^6 g'^2+70f'^2 g'^6+g'^8) \\
B_{8} &=&\frac{1}{4}c^{2}r^{8}f'g'(f'^6+7f'^4 g'^2+7f'^2 g'^4+g'^6)
\end{eqnarray*}
Since $\alpha (u)=(f(u),g(u))$ is not a constant planar curve, we
parametrize it by the arc-length, that is, $(f(u),g(u))=(x(\phi
(u),y(\phi (u)),$ where
$$f'(u)=\phi ^{\prime }(u)\cosh (\phi (u)),\hspace*{.5cm}g'(u)=\phi'(u)\sinh (\phi (u)),\hspace*{.5cm}\phi'^2=f'^2-g'^2.$$
With this change of variable, the functions $A_{8}$ and $B_{8}$
write now as:
$$A_{8} =-\frac{1}{32}c^{2}r^{8}\phi ^{\prime 8}\cosh (8\phi (u)),\hspace*{1cm}
B_{8} =\frac{1}{32}c^{2}r^{8}\phi ^{\prime 8}\sinh (8\phi (u)).$$
As $c\neq 0$ and $r>0$, we conclude that $\phi ^{\prime }=0$ on some
interval, that is, $\alpha $ is a constant curve, obtaining a
contradiction. This finishes the Theorem for the case that the foliation planes are timelike.
%%%%%%%%%%%%%%%%%%%%%%%%%%%%%%%%%%%%%%%%%%%%%

\subsection{The planes are lightlike }
%%%%%%%%%%%%%%%%%%%%%%%%%%%%%%%%%%%%%%%%%%%%%

After a motion in $\e_1^{3},$ we  parametrize the
surface by
\begin{equation}\label{para3}
\hbox{\bf
X}(u,v)=(f(u),g(u)+u,g(u)-u)+(v,r(u)\frac{v^{2}}{2},r(u)\frac{v^{2}}{2}),
\end{equation}
where $r>0, f$ and $g$ are smooth functions. In such case, $M$ is rotational if  $f$ is a constant function.

%%%%%%%%%%%%%%%%%%%%%%%%%%%%%%%%%%%%%
\subsubsection{Case $c=0$.}
%%%%%%%%%%%%%%%%%%%%%%%%%%%%%%%%%%%%%%

We compute (\ref{w1}) and we take $4b^2=1$ again. With our parametrization, and we
obtain
\begin{equation}\label{e-l}
\sum_{j=0}^{6}A_{j}(u)v^{n}=0,
\end{equation}
for some functions $A_{j}.$ As a consequence, all coefficients $A_{j}$ vanish. Then
$$A_{6}=-2a^2(2r^{2}-r')(-4rr'+r'')^{2}.$$
\begin{enumerate}
\item  If $2r^{2}-r^{\prime }=0$ then $r$ is given by
$$r(u)=\frac{1}{-2u-\lambda },\ \ \lambda \in\r.$$
Now
$$A_{3}=\frac{16a^2 f'(-4f'+(2u+\lambda )f'')^{2}}{(2u+\lambda )^{5}}.$$
From $A_{3}=0$ we have
\begin{enumerate}
\item  If $f'=0$ then $f$ is constant and $M$ is rotational.
\item If  $-4f'+(2u+\lambda )f''=0$ then $f''=\frac{4f'}{(2u+\lambda )}$. Putting it into $A_{2}$, it gives
$A_{2}=-256f'^2$, which implies that $f'=0$ and $M$  is rotational again.
\end{enumerate}
\item  Assume $-4rr'+r''=0$. The coefficient $A_4$ gives $a^2(rf''+2r'f')^2+2r'^2(2r^2-r')=0$.
A first integral of $-4rr'+r''=0$ is $2r^2-r'=k$, for some constant $k\not=0$. Then $A_4=0$ writes
$$a^2(rf''+2r'f')^2+2r'^2 k=0.$$
\begin{enumerate}
\item If $k>0$, then $r$ is constant and $f''=0$. In particular, $f(u)=\lambda u+\mu$.  Now $A_2=0$ implies $-16a^2 r^2(\lambda^2 r+4r g'+g'')^2=0$. Solving for $g$, we obtain $g(u)=-\lambda^2u/4-e^{-4ru c_1}/(4r)+c_2$, $c_1,c_2\in\r$. Hence, (\ref{e-l}) writes
$-256c_1r^4 e^{-8ru}=0$: contradiction.
\item Assume $k=-\lambda<0$. Then $A_4=0$ implies
$$rf''+2r'f'=\pm\sqrt{2\lambda/a^2}r'.$$
Here we obtain $f''$, which it is substituted in $A_3$ to obtain $g''$ in terms of $f'$ and $g'$. Substituting into $A_2$, we get
that $A_1=0$ is equivalent to
$(2r^2+k)/(2-2r^2)=0$. Thus, the only possibility is that $r$ is a constant function. But then $r'=2r^2+\lambda$ gets a contradiction.
\end{enumerate}

\end{enumerate}
%%%%%%%%%%%%%%%%%%%%%%%%%%%%%%%%%%%%%%%%%%%%%
\subsubsection{Case c$\neq 0$}
%%%%%%%%%%%%%%%%%%%%%%%%%%%%%%%%%%%%%%%%%%%%%

If we compute the Weingarten relation (\ref{w1}) with our parametrization,  we
obtain
\begin{equation}
\sum_{j=0}^{8}B_{j}(u)v^{n}=0,
\end{equation}
for some functions $B_{j}$. As a consequence, all coefficients $B_{j}$ vanish. The leader coefficient $B_8$ is
$$B_{8}=-64c^{2}(-2r^{2}+r')^{4}.$$
Thus $-2r^{2}+r'=0$ and $r$ is given by
$$r(u)=\frac{1}{-2u-\mu },\ \mu \in \r.$$
Now
$$B_{3}=\frac{1024c^{2}f^{\prime 4}}{(2u+\mu )^{5}}.$$
From $B_{3}=0$ we have $f'=0$ and thus $M$ is rotational.

\begin{remark}\label{r-1} The non-rotational spacelike surfaces in $\e_1^3$ with $H=0$ and $K=0$ are determined by the computation of (\ref{p0}). In the case $H=0$ and if the planes of the foliations are spacelike or timelike, the functions $f$, $g$ and $r$ in the parametrizations (\ref{para1}) and (\ref{para2}) satisfy
$$f'=\lambda r^2,\hspace*{.5cm}g'=\mu r^2,\hspace*{.5cm} 1-(\epsilon\lambda^2+\mu^2)r^4-r'^2+rr''=0,$$
with $\lambda, \mu\in\r$ and $\epsilon=1$ or $\epsilon=-1$ depending if the planes are spacelike or timelike, respectively.
The solutions are given in terms of elliptic equations. If the planes are lightlike, then $H=0$ means, up constants, $r=\tan(2u)$ and
$$f=\lambda(u+\frac12\cot(2u))\hspace*{.5cm}g=\frac{1}{32}\Big(4(\mu-3\lambda^2)u-4\lambda^2\cot(2u)
-(\lambda^2-4\mu)\sin(4u)\Big),\ \ \lambda,\mu\in\r.$$

If $K=0$, then the functions satisfy $f''=g''=r''=0$ when the foliation planes are spacelike or timelike, and
$f''=g''=0$ and $r=\lambda/(u+\mu)$ if the planes are lightlike.
\end{remark}

%%%%%%%%%%%%%%%%%%%%%%%%%%%%%%%%%%%%%%%%%%%%%%%%%%%

\section{Proof of Theorem \ref{t2}}

%%%%%%%%%%%%%%%%%%%%%%%%%%%%%%%%%%%%%%%%%%%%%%%%%%%%

Let $M$ be a linear Weingarten spacelike surface foliated by a uniparametric family of circles. Consider a real interval $I\subset\r$ and $u\in I$ the parameter of each plane of the foliation that defines $M$. Let $\G(u)$ be a smooth unit vector field orthogonal to each $u$-plane. Assume that the $u$-planes are not parallel and we will conclude that $M$ is a pseudohyperbolic surface. Then $\G'(u)\not=0$ in some real interval. Without loss of generality, we assume that in that interval, the planes containing the circles of $M$ have the same causal character. Consider an integral curve $\Gamma$ of the vector field $\G$. Then $\Gamma$ is not a straight-line. This allows to define a Frenet frame of $\Gamma$ $\{\T,\N,\B\}$. We distinguish three cases according to the causal character of the foliation planes.

\subsection{The planes are spacelike}

Let $\{\E_1(u),\E_2(u)\}$ be an orthonormal basis in each $u$-plane. Then $M$ parametrizes as
$$\X(u,v)=\C(u)+r(u)(\cos(v)\E_1(u)+\sin(v)\E_2(u))$$
where $r(u)>0$ and $\C(u)$ are differentiable functions on $u$.
Then $\T=\G$ is the unit tangent vector to $\Gamma$ and the Frenet equations are
\begin{eqnarray*}
\T'&=&\kappa \N\\
\N'&=&\kappa\T+\sigma \B\\
\B'&=&-\sigma\N
\end{eqnarray*}
A change of coordinates allows to write $M$ as
$$\X(u,v)=\C(u)+r(u)(\cos(v)\N(u)+\sin(v)\B(u)).$$
Set $\C'=\alpha\T+\beta\N+\gamma\B$, where $\alpha,\beta$ y $\gamma$ are smooth functions on $u$. Here $\T$ is a timelike unit vector and
$\N$ and $\B$ are spacelike unit vectors. Note that $\kappa\not=0$ since $\Gamma$ is not a straight-line.

By using $\C'$ and the Frenet equations, the expression (\ref{w1}) is a  trigonometric polynomial on $\cos(jv)$ and $\sin(jv)$. Exactly, there exist smooth functions on $u$, namely $A_j, B_j$, such that Equation (\ref{w1}) writes as
$$\sum_{j=0}^{8}A_{j}(u)\cos {(jv)}+B_{j}(u)\sin {(jv)}=0.$$

\subsubsection{Case $c=0$ in the relation $aH+bK=c$.}

Without loss of generality, we assume that $4b^2=1$.
 The coefficient $B_8$ implies
$$\beta\gamma\Big(2 a^2(3\beta^4-10\beta^2\gamma^2+3\gamma^4)+\kappa^2(1+12a^2 r^2)(\gamma^2-\beta^2)+r^2\kappa^4(1+6a^2r^2)\Big)=0.$$
We discuss three cases.

\begin{enumerate}
\item Case $\beta=0$ in a sub-interval of $I$. Then $A_8=0$ writes as
$$(\gamma^2+r^2\kappa^2)^2\Big(4a^2\gamma^2+(1+4a^2 r^2)\kappa^2\Big)=0.$$
Since $r\kappa\not=0$, this gives a contradiction.
\item Case $\gamma=0$ in a sub-interval of $I$. Equation $A_8=0$ writes as
$$(\beta^2-r^2\kappa^2)^2\Big((1+4a^2r^2)\kappa^2-4a^2\beta^2\Big)=0.$$
If $\beta^2=r^2\kappa^2$, it follows that
$$A_6=-\frac{9}{32}\kappa^6 r^{10}(\alpha-r')^2.$$
Then $A_6=0$ yields $\alpha=r'$. A computation of $W$ gives $W=0$: contradiction. As a consequence,  we assume $4\beta^2=(1+4a^2r^2)\kappa^2$. The computation of the coefficient $A_7$ leads to
$\alpha^2(1+4a^2r^2)=4a^2r^2r'^2$. From the expression of the center curve $\C$, we have
$$\C'=\frac{rr'}{\sqrt{\frac{1}{4a^2}+r^2}}\ \T+\kappa \sqrt{\frac{1}{4a^2}+r^2}\ \N=(\sqrt{\frac{1}{4a^2}+r^2}\ \T)'.$$
 In particular, there exists $\C_0\in\e_1^3$ such
that $$\C=\C_0+\sqrt{\frac{1}{4a^2}+r^2}\ \T.$$
The parametrization $\X$ of the surface is now
$$\X(u,v)=\C_{0}+\sqrt{\frac{1}{4a^2}+r^{2}}\ \T+r(\cos (v)\N+\sin (v)\B).$$
Then
$$\langle\X-\C_0,\X-\C_0\rangle=\frac{1}{4a^2}.$$
This means that the surface is a pseudohyperbolic surface.
\item Case $\beta\gamma\not=0$.  From $B_{8}=0$  we can
calculate $\beta ^{2}$:
$$\beta ^{2}=\frac{1}{12a^2}\Big(20a^2\gamma ^{2}+(1+12a^2r^{2})\kappa ^{2}\pm A\Big),$$
where $A=\sqrt{256a^{2}\gamma ^{4}+16a^2\gamma ^{2}\kappa
^{2}+192a^{2}r^{2}\gamma ^{2}\kappa ^{2}+\kappa ^{4}}.$ We consider the sign
'+' in the value of $\beta ^{2}$ (similarly with the choice '-'). Let us put it
into $A_{8}$ and taking in account that $\kappa \neq 0$, we obtain the following identity
\begin{eqnarray*}
& &26624\gamma ^{6}+\kappa ^{6}+1536a^{2}\gamma ^{4}\kappa
^{2}(1+2a^2r^{2})+ 72a^2\gamma ^{2}\kappa ^{4}\\
&=&-(1792a^{2}\kappa ^{4}+\kappa
^{4}+64a^2\gamma ^{2}\kappa ^{2}(1+12a^2r^{2}))A.
\end{eqnarray*}
Squaring both sides and after some manipulations, we obtain
$$(\gamma ^{2}+\kappa ^{2}r^{2})\Big((16a^2\gamma ^{2}+\kappa
^{2})^{2}+256a^{2}r^{2}\gamma ^{2}\kappa ^{2}\Big)=0.$$
This would be imply $\kappa r=0$: contradiction.
\end{enumerate}

\subsubsection{Case $c\not=0$ in the relation $aH+bK=c$.}

Without loss of generality, we shall assume that $c=1$. Also, we discard the cases $b=0$ or $a=0$, corresponding to the known situations of non-zero constant mean curvature or constant Gauss curvature. The computation of
the coefficients $A_{8}$ and $B_{8}$ gives
$$A_{8}=-\frac{1}{32}r^{8}x_{1},\hspace*{1cm}B_{8}=\frac{1}{16}\beta \gamma
r^{8}x_{2},$$
where
\begin{eqnarray*}
x_{1} &=&\beta ^{8}-(28\gamma ^{2}+\kappa ^{2}(a^2+2b+4r^{2}))\beta ^{6}
\label{6} \\
&&+(70\gamma ^{4}+15\gamma ^{2}\kappa ^{2}(^{2}(a^2+2b+4r^{2})+\kappa
^{4}(b^{2}+3(a^2+2b)r^{2}+6r^{4}))\beta ^{4}   \\
&&+(-28\gamma ^{6}-15\gamma ^{4}\kappa ^{2}(a^2+2b+4r^{2})-\kappa
^{6}r^{2}(2b^{2}+3(a^2+2b)r^{2}+4r^{4})\\
& &-6\gamma ^{2}\kappa
^{4}(b^{2}+3(a^2+2b)r^{2}+6r^{4}))\beta ^{2}   \\
&&+(\gamma ^{2}+r^{2}\kappa ^{2})^{2}(\gamma ^{4}+\gamma ^{2}\kappa
^{2}(a^2+2b+2r^{2})+\kappa ^{4}(b^{2}+(a^2+2b)r^{2}+r^{4})).
\end{eqnarray*}
\begin{eqnarray*}
x_{2} &=&-4\beta ^{6}+(28\gamma ^{2}+3\kappa ^{2}(m+2b+4r^{2}))\beta ^{4}
\label{7} \\
&&-2(14\gamma ^{4}+5\gamma ^{2}\kappa ^{2}(a^2+2b+4r^{2})+\kappa
^{4}(b^{2}+3(a^2+2b)r^{2}+6r^{4})\beta ^{2} \\
&&+(\gamma ^{2}+r^{2}\kappa ^{2})(4\gamma ^{4}+\gamma ^{2}\kappa
^{2}(3a^2+6b+8r^{2})+\kappa ^{4}(2b^{2}+3(a^2+2b)r^{2}+4r^{4}).
\end{eqnarray*}
From $B_8=0$, we discuss three cases:
\begin{enumerate}
\item Case $\gamma=0$ in some sub-interval of $I$. Then $A_8=0$ is
$$(\beta^2-r^2\kappa^2)^2\Big(\beta^4-(a^2+2b+2r^2)\beta^2\kappa^2+(b^2+(a^2+2b)r^2+r^4)\kappa^4\Big)=0.$$
\begin{enumerate}
\item Suppose $\beta^2=\kappa^2r^2$. Without loss of generality we
assume that $\beta =\kappa r$. Now $A_{6}=-\frac{9}{8}b^{2}\kappa
^{6}r^{10}(\alpha -r^{\prime })^{2}$. Then $\alpha =r^{\prime }$ and the computation of $W$ gives $W=0:$ contradiction.
\item Then $\beta^4-(a^2+2b+2r^2)\beta^2\kappa^2+(b^2+(a^2+2b)r^2+r^4)\kappa^4=0$. If we look this expression as a polynomial on $\beta^2$, from the discriminant we conclude that $a^2+4b\geq 0$. If $a^2+4b=0$, $\beta^2=(a^2+4r^2)\kappa^2/4$. Then $B_5=0$ gives
 $$ B_{5}=\frac{1}{128}a^{4}\kappa ^{5}r^{7}\sigma \sqrt{a^{2}+4r^{2}}(\alpha
\sqrt{a^{2}+4r^{2}}-2rr^{\prime })^{2}=0. $$
If $\sigma=0$, $A_5=0$ implies $\alpha\sqrt{a^{2}+4r^{2}}-2rr^{\prime }=0$ again. Therefore, in both cases, and from the value of $\alpha$, we can write
$$\C'=\left( \frac{\sqrt{a^{2}+4r^{2}}}{2}\ \T\right)'$$
and so, there exists $\C_0\in\e_1^3$ such that
$$\C=\C_0+\frac{\sqrt{a^{2}+4r^{2}}}{2}\ \T.$$
As a consequence, we have again
$$\X(u,v)=\C_0+\frac12\sqrt{a^2+4r(u)^2}\ \T+r(u) (\cos(v)\N+\sin(v)\B,$$
for some $\C_0\in\e_1^3$. Therefore $\langle\X-\C_0,\X-\C_0\rangle=-a^2/4$,
and the surface is a pseudohyperbolic surface.

Assume then $a^2+4b>0$. The coefficient $A_{7}$ is
$$A_{7}=\frac{1}{64}aAB\kappa ^{5}r^{9}(\alpha \kappa ^{2}-\kappa \beta
^{\prime }+\kappa ^{\prime }\beta )=0,$$
with
$$A=2b+a(a+\sqrt{a^{2}+4b}),\hspace*{1cm} B=a^{3}+4ab+(a^{2}+2b)\sqrt{a^{2}+4b}.$$
Then  number $A$ does not vanish and $B=0$ holds only if $a^{2}+4b=0$. From $A_7=0$ we conclude that
$\alpha \kappa ^{2}-\kappa \beta
^{\prime }+\kappa ^{\prime }\beta =0$, that is,
$$\alpha =\left( \frac{\beta }{\kappa }\right) ^{\prime },$$
which implies $\C=\C_0+\beta/\kappa\T$ for some $\C_0\in\e_1^3$.
The derivative of the curve $\C$ is
$$\C^{\prime }=\left( \frac{\beta }{\kappa }\right) ^{\prime }\T+\beta \N=\left(
\frac{\beta }{\kappa }\T\right) ^{\prime }.$$
The expression of $\X(u,v)$ is
$$\X(u,v)=\C_{0}+\frac{\beta }{\kappa }\T+r(\cos(v) \N+\sin(v)\B).$$
 Using the value of $\beta^{2}$, we have,
$$
\langle\X-\C_0,\X-\C_0\rangle=-\frac{\beta^{2}}{\kappa ^{2}}+r^{2}=-\left( \frac{a^{2}}{2}+b+\frac{a}{2}\sqrt{a^{2}+4b}\right).$$
This means that the surface is a  pseudohyperbolic surface.
\end{enumerate}
\item Case $\beta=0$ in some sub-interval of $I$. Then
$$A_8=-\frac{1}{32}r^8 (\gamma^2+\kappa^2 r^2)y_1,\hspace*{1cm}A_7=-\frac{1}{16}\alpha\kappa r^9(\gamma^2+\kappa^2 r^2)z_1,$$
where
\begin{eqnarray*}y_1&=&\gamma^4+(a^2+2b+2r^2)\kappa^2\gamma^2+(b^2+(a^2+2b)r^2+r^4)\kappa^4\\
z_1&=&8\gamma^4+(7(a^2+2b)+16r^2)\kappa^2\gamma^2+(6b^2+7(a^2+2b)r^2+8r^4)\kappa^4.
\end{eqnarray*}
Assume $\alpha\not=0$. From $y_1=0$, we obtain $\gamma^2$, which it is substituted into $z_1=0$, obtaining
$a\sqrt{a^2+4b}=\pm (a^2+2b)$. Then $a^2(a^2+4b)=(a^2+2b)^2$, which implies $b=0$: contradiction. Therefore,
$\alpha=0$. From $y_1=0$,
$$\gamma^4+(a^2+2b+2r^2)\kappa^2\gamma^2+(b^2+(a^2+2b)r^2+r^4)\kappa^4=0.$$
Then
\begin{equation}\label{gamma}
\gamma^2=\frac12\Big(\pm a\sqrt{a^2+4b}-(a^2+2b+2r^2)\Big)\kappa^2.
\end{equation}
We prove that the quantity in the parenthesis is non-positive, that is,
$\pm a\sqrt{a^2+4b}-(a^2+2b+2r^2)\leq 0$. Since this function on $r$ is decreasing on $a$, we show that (taking $r=0$)
$\pm a\sqrt{a^2+4b}-(a^2+2b)\leq 0$. Depending on the sign of $a$, we have two possibilities. If $a>0$, the inequality
$\pm a\sqrt{a^2+4b}\leq a^2+2b$ is trivial. If $a<0$, the inequality is trivial if $a^2+2b\geq 0$. The only case to consider is
$\pm a\sqrt{a^2+4b}\leq a^2+2b<0$ ($\Rightarrow b<0$). But $a^2+2b<0$ and $a^2+4b\geq 0$ is not compatible.
As a consequence of this reasoning, we conclude from (\ref{gamma}) that $\gamma=0$. This case was  studied in the above subsection.

\item Case $\beta\gamma\not=0$. In this last case, the computations become very complicated and difficult. For this reason, we only give the proof outline  and we omit the details.
    Let $x=\beta^2$, $y=\gamma^2$. From $x_1=0$, we obtain the value of $a^2+2b$, which is substituted into $x_2=0$, obtaining
    $$\Big((x + y)^2+2(y-x)r^2\kappa^2+r^4\kappa^4\Big)^2\Big((x+y)^2+2(y-x)r^2\kappa^2-b^2\kappa^4+r^4\kappa^4\Big)=0.$$
    If we see $(x + y)^2+2(y-x)r^2\kappa^2+r^4\kappa^4=0$ as polynomical equation on $r^2\kappa^2$, we find that the discriminant is negative, and so, this case is impossible. Thus
    \begin{equation}\label{q}
    (x+y)^2+2(y-x)r^2\kappa^2-b^2\kappa^4+r^4\kappa^4=0.
    \end{equation}
Then
$y=-x-r^2\kappa^2+\kappa\sqrt{4xr^2+b^2\kappa^2}$. Putting into $x_1=0$, we conclude
$$16x^2-8x(a^2+2b+2r^2)\kappa^2+(a^4+4a^2b)\kappa^4=0$$
or
$$256x^4-512x^3r^2\kappa^2-
128x^2(b^2-2r^4)\kappa^4+64b^2xr^2\kappa^6+3b^4\kappa^8=0.$$
We analyse the first possibility (the second one is analogous).
If $16x^2-8x(a^2+2b+2r^2)\kappa^2+(a^4+4a^2b)\kappa^4=0$, then
\begin{eqnarray}
\beta^2&=&\frac{\kappa^2}{4}\Big(a^2+2b+2r^2\pm 2Q\Big)\label{bg1}\\
\gamma^2&=&\frac{\kappa^2}{4}
\Big(4\sqrt{b^2+r^2(a^2+2b+2r^2 \pm 2Q)}-(a^2+2b+6r^2) \mp 2Q\Big),\label{bg2}
\end{eqnarray}
where
$$Q=\sqrt{a^2r^2+(b+r^2)^2)}.$$
With these values obtained for $\beta^2$ and $\gamma^2$, we place them into $x_1=0$, obtaining
an equation that depends only on the function $r=r(u)$. Exactly, it is  a rational expression on  $r$ and $\sqrt{b^2+(a^2+2b)r^2+r^4}$
$${\cal P}(r,\sqrt{b^2+(a^2+2b)r^2+r^4})=0.$$
In particular, $r(u)$ is a constant function.

We do the change
    $$p=x-y,\hspace*{1cm}q=(x-y)^2-4xy$$
    that is, $x=(p+\sqrt{2p^2-q})/2$ and $y=(-p+\sqrt{2p^2-q})/2$. Equation (\ref{q}) writes as
    \begin{equation}\label{pq}
    2p^2-q-2pr^2\kappa^2+(r^4-b^2)\kappa^4=0\Rightarrow q=2p^2+2pr^2\kappa^2+(r^4-b^2)\kappa^4.
    \end{equation}
We calculate $q$ by other way. From $x_1=0$, we get a value of $q$, which it is substituted into $x_2=0$, obtaining
$$(p-r^2\kappa^2)\Big(2p-(a^2+2(b+r^2))\kappa^2\Big)\Big(4p^2-2p(a^2+2b+4r^4)\kappa^2+(b^2+2a^2r^2+4br^2+4r^4)
\kappa^4\Big)=0.$$
Hence and together (\ref{pq}) we obtain different values for $p$ and $q$.
    On the other hand, the values obtained for $\beta$ and $\gamma$ in (\ref{bg1}) and (\ref{bg2}) allow to get a pair of values for $p$ and $q$, which they must be equal. For each pair of these values, we obtain
       different values for $r$, which are substituted into the coefficients $A_i$ and $B_i$. The computation of the coefficients $A_7$ and $B_7$ gives $\alpha=0$ and using $A_5$ and $B_5$, we get  $\sigma=0$. Finally the coefficients    $A_4$ and $B_4$ give $\kappa=0$, obtaining a contradiction.

\end{enumerate}
\subsection{The planes are timelike}
This situation is similar that the above case of spacelike planes and we omit the details.

\subsection{The planes are lightlike}

The surface $M$ can be locally written as
$$\X(u,v)=\C(u)+v\N(u)+r(u) v^2 \T(u),$$
where $r(u)>0$, and $\T$ and $\N$ are the tangent vector and normal vector of $\Gamma$ respectively. Since the planes are lightlike,
$\langle\T,\T\rangle=0$ and $\langle\N,\N\rangle=1$. The Frenet frame for $\Gamma$ is
$\{\T,\N,\B\}$, where $\B$ is the unique lightlike vector orthogonal to $\N$ such that $\langle\T,\B\rangle=1$ and
$\mbox{det}(\T,\N,\B)=1$.
The Frenet equations are
\begin{eqnarray*}
\T'&=&\kappa \N\\
\N'&=&\sigma\T-\kappa \B\\
\B'&=&-\sigma\N
\end{eqnarray*}
Again, we put $\C'=\alpha\T+\beta\N+\gamma\B$. By using $\C'$ and the Frenet equations, the expression (\ref{w1}) is a  trigonometric polynomial on $v$ such as $\sum_{j=0}^{n}A_{j}(u)v^j=0$ with $n=11$ if $c=0$ and $n=12$ if $c\not=0$.

\subsubsection{Case $c=0$ in the relation $aH+bK=c$.}
Without loss of generality, we assume that $b=1/2$. Then
$$A_{11}=98a^2r^2\kappa^5(2r^2\gamma-r')^2.$$
With this value of $r'$, we obtain
$$A_8=-64r^2\kappa^5(\sigma-2r\beta)^2(-4a^2r\beta+r^2\kappa+2a^2\sigma)=0.$$
If $\sigma=2r\beta$, then $A_6=-100r^4\alpha^2\kappa^6$. This yields $\alpha=0$ and $W=0$: contradiction. Thus
$2a^2 \sigma=4a^2r\beta-r^2\kappa$. Now
$A_7=0$ gives $2a^2\alpha=r^2\gamma$. Then
$$\C'=\frac{r'}{4a^2}\T+\beta\N+\frac{r'}{2r^2}\B=\Big(\frac{r}{4a^2}\T-\frac{1}{2r}\B\Big)'.$$
Therefore, there exists $\C_0\in\e_1^3$ such that
$$\X(u,v)=\C_0+\Big(\frac{r}{4a^2}\T-\frac{1}{2r}\B\Big)+v\N(u)+r v^2 \T(u).$$
In particular,
$$\langle\X(u,v)-\C_0,\X(u,v)-\C_0\rangle=-\frac{1}{4a^2},$$
which shows that the surface is a pseudohyperbolic surface.

\subsubsection{Case $c\not=0$ in the relation $aH+bK=c$.}
We assume that $c=1$. Then
$$A_{12}=-64\kappa^4(r'-2r^2\gamma)^4.$$
As above, $A_8=0$ gives two possibilities about the value of $\sigma$. In the first case,
$\sigma=2r\beta$ and $A_6=0$ yields $\alpha=0$. This implies $W=0$: contradiction. The other case for $\sigma$ is
$$\sigma^2+2r\sigma (-2\beta+(a^2+2b)r\kappa)+4(\beta^2-(a^2+2b)r\beta\kappa+b^2r^2\kappa^2)r^4=0.$$
Then
$$\sigma=2r\beta-a^2r^2\kappa-2br^2\kappa\pm ar^2\kappa\sqrt{a^2+4b}.$$
In particular, $a^2+4b\geq 0$. We assume the choice '+' in the above identity (similar for '-'). From Equation $A_7=0$, we obtain
$$\alpha=\frac12(a^2+2b-a\sqrt{a^2+4b})r'.$$
As in the case $c=0$, we conclude the existence of $\C_0\in\e_1^3$ such that
$$\C=\C_0+\frac12(a^2+2b-a\sqrt{a^2+4b})r\T-\frac{1}{2r}\B.$$
Now
$$\X(u,v)=\C_0+\frac12(a^2+2b-a\sqrt{a^2+4b})r\T-\frac{1}{2r}\B+v\N(u)+r v^2 \T(u),$$
and
$$\langle\X(u,v)-\C_0,\X(u,v)-\C_0\rangle=-\frac12(a^2+2b-a\sqrt{a^2+4b}),$$
showing that $M$ is a pseudohyperbolic surface again.
%%%%%%%%%%%%%%%%%%%%%%%%%%%%%%%%%%%%%%

\end{document}